\documentclass[12pt]{article}

\usepackage{amsmath,amsthm,amsfonts,amssymb,latexsym}
\usepackage{amsmath}

\usepackage{a4wide}
\newtheorem{thm}{Theorem}[section]

\newtheorem{lem}[thm]{Lemma}

\newtheorem{theorem}{Theorem}[section]

\newtheorem{proposition}[theorem]{Proposition}

\begin{document}

%-----------------------------------------------------------------------------------

\title{Cohomology of $\mathfrak {osp}(1|2)$ acting
on linear differential operators on the supercircle $S^{1|1}$ }

%\label{firstpage}

\author{Imed Basdouri\and Mabrouk Ben Ammar\thanks{D\'epartement de Math\'ematiques,
Facult\'e des Sciences de Sfax, BP 802, 3038 Sfax, Tunisie.
E.mails:basdourimed@yahoo.fr, mabrouk.benammar@fss.rnu.tn}}

%--------------------------------------------------------------------
\maketitle
% ----------------------------------------------------------------

\begin{abstract}

We compute the first cohomology spaces
$H^1\left(\mathfrak{osp}(1|2);\mathfrak{D}_{\lambda,\mu}\right)$
($\lambda,\,\mu\in\mathbb{R}$) of the Lie superalgebra
$\mathfrak{osp}(1|2)$ with coefficients in the superspace
$\mathfrak{D}_{\lambda,\mu}$ of linear differential operators
acting on weighted densities on the supercircle $S^{1|1}$. The
structure of these spaces was conjectured in \cite{gmo}. In fact,
we prove here that the situation is a little bit more complicated.
(To appear in LMP.)
\end{abstract}

\maketitle {\bf Mathematics Subject Classification} (2000). 53D55


%\thispagestyle{empty}

%----------------------------------------------------------------------------------

\section{Introduction}

%-----------------------------------------------------------------------------------

The space of weighted densities with weight $\lambda$ (or $\lambda$-densities) on
$S^1$, denoted by:
\begin{equation*}
{\mathcal F}_\lambda=\left\{ f(dx)^{\lambda}, ~f\in C^\infty(S^1)\right\},\quad (\lambda\in\mathbb{R}),
\end{equation*}
is the space of sections of the line bundle $(T^*S^1)^{\otimes^\lambda}.$
Let ${\rm Vect}(S^1)$ be the Lie algebra of all vector fields $F\frac{d}{dx}$
on $S^1$, ($F\in C^\infty(S^1)$). With the {\it Lie derivative}, ${\mathcal F}_\lambda$ is a
${\rm Vect}(S^1)$-module. Alternatively, the ${\rm Vect}(S^1)$ action can be written as follows:
\begin{equation}\label{Lie1}
L_{F\frac{d}{dx}}^\lambda(f(dx)^{\lambda})=(Ff'+\lambda fF')(dx)^{\lambda},
\end{equation}
where $f'$, $F'$ are $\frac{df}{dx}$, $\frac{dF}{dx}$.

Let $A$ be a differential operator on $S^1$. We see $A$ as the linear mapping
$f(dx)^\lambda\mapsto(Af)(dx)^\mu$ from $\mathcal{F}_\lambda$ to
$\mathcal{F}_\mu$ ($\lambda$, $\mu$ in $\mathbb R$).
Thus the space of differential operators is a ${\rm Vect}(S^1)$ module,
denoted $\mathcal D_{\lambda,\mu}$. The ${\rm Vect}(S^1)$ action is:
\begin{equation}\label{Lieder2}
L_X^{\lambda,\mu}(A)=L_X^\mu\circ A-A\circ L_X^\lambda.
\end{equation}

If we restrict ourselves to the Lie subalgebra of $\mathrm{Vect}(S^1)$
generated by $\left\{\frac{d}{dx},\,x\frac{d}{dx},\,x^2\frac{d}{dx}\right\}$,
isomorphic to $\mathfrak{sl}(2)$, we get a family of infinite dimensional
$\mathfrak{ sl}(2)$ modules, still denoted ${\mathcal D}_{\lambda,\mu}$.

P. Lecomte, in \cite{lec}, found the cohomology spaces
$\mathrm{H}^1\left(\mathfrak{sl}(2);{\mathcal D}_{\lambda,\mu}\right)$ and
$\mathrm{H}^2\left(\mathfrak{sl}(2);{\mathcal D}_{\lambda,\mu}\right)$.
These spaces appear naturally in the problem of describing the deformations of
the $\mathfrak{sl}(2)$-module ${\mathcal D}$ of the differential operators acting on
${\mathcal S}^n=\bigoplus_{k=-n}^n{\mathcal F}_{\frac{1+k}{2}}$. More precisely,
the first cohomology space $\mathrm{H}^1\left(\mathfrak{sl}(2);V\right)$ classifies
the infinitesimal deformations of a $\mathfrak{sl}(2)$ module $V$ and the obstructions
to integrability of a given infinitesimal deformation of $V$ are elements of
$\mathrm{H}^2\left(\mathfrak{sl}(2);V\right)$. Thus, for instance, the infinitesimal
deformations of the $\mathfrak{sl}(2)$ module ${\mathcal D}$  are
classified by:
\begin{equation*}
\mathrm{H}^1\left(\mathfrak{sl}(2);{\mathcal D}\right)=\oplus_{k=0}^n\mathrm{H}^1
\left(\mathfrak{sl}(2);{\mathcal D}_{{\frac{1-k}{2}},{\frac{1+k}{2}}}\right)\oplus
\oplus_{k=-n}^n\mathrm{H}^1\left(\mathfrak{sl}(2);{\mathcal D}_{{\frac{1+k}{2}},{\frac{1+k}{2}}}\right).
\end{equation*}

In this paper we are interested to the study of the corresponding super structures.
More precisely, we consider here the superspace $S^{1|1}$ equipped with its standard
{\it contact structure 1-form} $\alpha$, and introduce the superspace $\mathfrak{F}_\lambda$
of $\lambda$-densities on the supercircle $S^{1|1}$.

Let $\mathcal{K}(1)$ be the Lie superalgebra of contact vector
fields, $\mathfrak{F}_\lambda$ is naturally a $\mathcal{K}(1)$-
module. For each $\lambda$, $\mu$ in $\mathbb{R}$, any differential
operator on $S^{1|1}$ becomes a linear mapping from
$\mathfrak{F}_\lambda$ to $\mathfrak{F}_\mu$, thus the space of
differential operators becomes a $\mathcal{K}(1)$-module denoted
$\mathfrak{D}_{\lambda,\mu}$.

To the symplectic Lie algebra $\mathfrak{sl}(2)$ corresponds the
ortosymplectic Lie superalgebra $\mathfrak{osp}(1|2)$ which is
naturally realized as a subalgebra of $\mathcal{K}(1)$. Restricting
our $\mathcal{K}(1)$-modules to $\mathfrak{osp}(1|2)$, we get
$\mathfrak{osp}(1|2)$-modules still denoted $\mathfrak{F}_\lambda$,
$\mathfrak{D}_{\lambda,\mu}$.

We compute here the first cohomology spaces
$\mathrm{H}^1\left(\mathfrak
{osp}(1|2);\mathfrak{D}_{\lambda,\mu}\right) $, ($\lambda$, $\mu$ in
$\mathbb{R}$), getting a result very close to the classical spaces
$\mathrm{H}^1\left(\mathfrak{sl}(2);{\mathcal
D}_{\lambda,\mu}\right)$. Especially, these spaces have the same
dimension. Moreover, we give explicit formulae for all the non
trivial 1-cocycles.

These spaces arise in the classification of infinitesimal
deformations of the $\mathfrak{osp}(1|2)$-module of the
differential operators acting on $ {\mathfrak
S}^n=\bigoplus_{k=1-n}^n\mathfrak{F}_{\frac{k}{2}}$. We hope to be
able to describe in the future all the deformations of this
module.

%----------------------------------------------------------------------------------------------------------------
\section{Definitions and Notations}
%----------------------------------------------------------------------------------------------------------------
%-----------------------------------------------------------------------------------------------------------------
\subsection{The Lie superalgebra of contact vector fields on ${S}^{1|1}$}
%------------------------------------------------------------------------------------------------------------------

We define the supercircle ${S}^{1|1}$ through its space of functions, $C^\infty({S}^{1|1})$.
A $C^\infty({S}^{1|1})$ has the form:
\begin{equation*}
F(x,\theta)=f_0(x)+\theta f_1(x),
\end{equation*}
where $x$ is the even variable and $\theta$ the odd variable: we have $\theta^2=0$.
Even elements in $C^\infty({S}^{1|1})$ are the functions
$F(x,\theta)=f_0(x)$, the functions $F(x,\theta)=\theta f_1(x)$ are odd elements.
Note $p(F)$ the parity of a homogeneous function $F$.

Let $\mathrm{Vect}(S^{1|1})$ be the superspace of vector fields on $S^{1|1}$:
\begin{equation*}
\mathrm{Vect}(S^{1|1})=\left\{F_0\partial_x+F_1\partial_\theta \quad F_i\in C^\infty(S^{1|1})\right\},
\end{equation*}
where $\partial_\theta$ and $\partial_x$ stand for $\frac{\partial}{\partial\theta}$ and
$\frac{\partial}{\partial_x}$. The vector fields $f(x)\partial_x$, and $\theta f(x)\partial_\theta$
are even, the vector fields $\theta f(x)\partial_x$, and $f(x)\partial_\theta$ are odd.
The superbracket of two vector fields is bilinear and defined for two homogeneous vector fields by:
$$
[X,Y]=X\circ Y-(-1)^{p(X)p(Y)}Y\circ X.
$$
Denote $\mathfrak{L}_X$ the Lie derivative of a vector field, acting
on the space of functions, forms, vector fields,\dots

The supercircle ${S}^{1|1}$ is equipped with the standard contact structure given by the
following even $1$-form:
\begin{equation*}
\alpha=dx+\theta d\theta.
\end{equation*}

We consider the Lie superalgebra $\mathcal{K}(1)$ of contact vector fields on $S^{1|1}$.
That is, $\mathcal{K}(1)$ is the superspace of conformal vector
fields on $S^{1|1}$ with respect to the $1$-form $\alpha$:
$$
\mathcal{K}(1)=\big\{X\in\mathrm{Vect}(S^{1|1})~|~\hbox{there exists}~F\in C^\infty({S}^{1|1})~
\hbox{such that}~\mathfrak{L}_X(\alpha)=F\alpha\big\}.
$$

Let us define the vector fields $\eta$ and $\overline{\eta}$ by
$$
\eta=\partial_\theta+\theta\partial_x,\quad \overline{\eta}=\partial_\theta-\theta\partial_x.
$$
Then any contact vector field on $S^{1|1}$ can be written in the following explicit form:
\begin{equation*}
X_F=F\partial_x+\frac{1}{2}\eta(F)(\partial_\theta-\theta\partial_x)=-F\overline{\eta}^2+
\frac{1}{2}{\eta}(F)\overline{\eta},\;\text{ where }\, F\in C^\infty(S^{1|1}).
\end{equation*}
Of course, $\mathcal{K}(1)$ is a subalgebra of
$\mathrm{Vect}(S^{1|1})$, and $\mathcal{K}(1)$ acts on
$C^{\infty}(S^{1|1})$ through:
\begin{equation}
\label{superaction3}
\mathfrak{L}_{X_F}(G)=FG'+\frac{1}{2}(-1)^{(p(F)+1)p(G)}\overline{\eta}(F)\cdot\overline{\eta}(G).
\end{equation}

Let us define the contact bracket on $C^\infty({S}^{1|1})$ as the bilinear mapping such that,
for a couple of homogenous functions $F$, $G$,
\begin{equation}
\{F,G\}=FG'-F'G+\frac{1}{2}(-1)^{p(F)+1}\overline{\eta}(F)\cdot\overline{\eta}(G),
\end{equation}

Then the bracket of $\mathcal{K}(1)$ can be written as:
\begin{equation*}
[X_F,\,X_G]=X_{\{F,\,G\}}.
\end{equation*}

%-----------------------------------------------------------------------------------------------
\subsection{The superalgebra $\mathfrak{osp}(1|2)$}
%----------------------------------------------------------------------------------------------

Recall the Lie algebra $\mathfrak{sl}(2)$ is isomorphic to the Lie subalgebra of ${\rm Vect}(S^1)$
generated by
\begin{equation*}\left\{\frac{d}{dx},\,x\frac{d}{dx},\,x^2\frac{d}{dx}\right\}.
\end{equation*}
Similarly, we now consider the orthosymplectic Lie superalgebra as a
subalgebra of $\mathcal{K}(1)$:
\begin{equation*}
\mathfrak{osp}(1|2)=\text{Span}(X_1,\,X_{x},\,X_{x^2},\,X_{x\theta},\,
X_{\theta}).
\end{equation*}
The space of even elements is isomorphic to $\mathfrak{sl}(2)$:
\begin{equation*}
(\mathfrak{osp}(1|2))_0=\text{Span}(X_1,\,X_{x},\,X_{x^2})=\mathfrak{sl}(2).
\end{equation*}
The space of odd elements is two dimensional:
\begin{equation*}
(\mathfrak{osp}(1|2))_1=\text{Span}(X_{x\theta},\,X_{\theta}).
\end{equation*}
The new commutation relations are
\begin{equation*}
\aligned
&[X_{x^2},X_\theta]=-X_{x\theta},~~&[X_x,X_\theta]=-\frac{1}{2} X_\theta,~~&[X_1,X_\theta]=0,\cr
&[X_{x^2},X_{x\theta}]=0,~~&[X_x,X_{x\theta}]=\frac{1}{2} X_{x\theta},~~&[X_1,X_{x\theta}]=X_\theta,\cr
&[X_{x\theta},X_{\theta}]=\frac{1}{2} X_{x}.&&
\endaligned
\end{equation*}

%----------------------------------------------------------------------------------------
\subsection{The space of weighted densities on $S^{1|1}$}
%----------------------------------------------------------------------------------------

In the super setting, by replacing $dx$ by the 1-form $\alpha$, we get analogous definition for
weighted densities {\sl i.e.} we define the space of $\lambda$-densities as
\begin{equation}
\mathfrak{F}_\lambda=\left\{\phi=F(x,\theta)\alpha^\lambda~~|~~F(x,\theta) \in C^\infty(S^{1|1})\right\}.
\end{equation}
As a vector space, $\mathfrak{F}_\lambda$ is isomorphic to $C^\infty(S^{1|1})$,
but the Lie derivative of the density $G\alpha^\lambda$ along the vector field $X_F$ in
$\mathcal{K}(1)$ is now:
\begin{equation}
\label{superaction}
\mathfrak{L}_{X_F}(G\alpha^\lambda)=\mathfrak{L}^{\lambda}_{X_F}(G)\alpha^\lambda,
\quad\text{with}~~\mathfrak{L}^{\lambda}_{X_F}(G)=\mathfrak{L}_{X_F}(G)+ \lambda F'G.
\end{equation}
Or, if we put $F=a(x)+b(x)\theta$, $G=g_0(x)+g_1(x)\theta$,
\begin{equation}
\label{deriv}
\mathfrak{L}^{\lambda}_{X_F}(G)=L^{\lambda}_{a\partial_x}(g_0)+\frac{1}{2}~bg_1+
\left(L^{\lambda+\frac{1}{2}}_{a\partial_x}(g_1)+\lambda
g_0b'+\frac{1}{2} g'_0 b\right)\theta.
\end{equation}

Especially, we have
\begin{equation*}
\left\{\begin{array}{lll}
\mathfrak{L}^\lambda_{X_a}(g_0)=L^\lambda_{a\partial_x}(g_0),&&
\mathfrak{L}^\lambda_{X_a}(g_1\theta)=\theta L^{\lambda+\frac{1}{2}}_{a\partial_x} (g_1),\cr
&\hbox{ ~~and~~}&\cr
\mathfrak{L}^{\lambda}_{X_{b\theta}}(g_0)=(\lambda g_0b'+\frac{1}{2} g'_0b)
\theta&&\mathfrak{L}^{\lambda}_{X_{b\theta}}(g_1\theta)=\frac{1}{2} bg_1.
\end{array}\right.\end{equation*}

Of course, for all $\lambda$, $\mathfrak{F}_\lambda$ is a
$\mathcal{K}(1)$-module:
\begin{equation*}
[\mathfrak{L}^{\lambda}_{X_F},\mathfrak{L}^{\lambda}_{X_G}] =
\mathfrak{L}^{\lambda}_{[{X_F},\,X_G]}.
\end{equation*}
We thus obtain a one-parameter family of $\mathcal{K}(1)$-modules on
$C^\infty(S^{1|1})$ still denoted by $\mathfrak{F}_\lambda$.

%------------------------------------------------------------------------------------------
\subsection{Differential Operators on Weighted Densities}

%----------------------------------------------------------------------------------
A differential operator on $S^{1|1}$ is an operator on $C^\infty(S^{1|1})$ of t
he following form:
$$
A= \sum_{i=0}^\ell\widetilde{a}_i (x, \theta)\partial_x^i +
\sum_{i=0}^\ell\widetilde{b}_i (x, \theta)\partial_x^i\partial_\theta.
$$
In \cite{gmo}, it is proved that any local operator $A$ on $S^{1|1}$ is in
fact a differential operator.

Of course, any differential operator defines a linear mapping from $\mathfrak{F}_\lambda$
to $\mathfrak{F}_\mu$ for any $\lambda$, $\mu\in\mathbb{R}$, thus the space of differential
operators becomes a family of $\mathcal{K}(1)$ and $\mathfrak{osp}(1|2)$ modules denoted
$\mathfrak{D}_{\lambda,\mu}$, for the natural action:
\begin{equation}\label{d-action}
\mathfrak{L}^{\lambda,\mu}_{X_F}(A)=\mathfrak{L}^\mu_{X_F}\circ A-
(-1)^{p(A)p(F)}A\circ \mathfrak{L}^\lambda_{X_F}.
\end{equation}

%-------------------------------------------------------------------------------------------------------------

%-------------------------------------------------------------------------------------------------------------

\section{The space $H^1(\mathfrak {osp}(1|2);\mathfrak{D}_{\lambda,\mu})$}

%----------------------------------------------------------------------------------------------------------------
%----------------------------------------------------------------------------------------------------------------

%--------------------------------------------------------------------------------------------------------------
\subsection{Lie superalgebra cohomology (see \cite{Fu})}

%--------------------------------------------------------------------------------------------------------------------

Let $\mathfrak{g}=\mathfrak{g}_0\oplus \mathfrak{g}_1$ be a Lie
superalgebra and $A=A_0\oplus A_1$ a $\mathfrak{g}$ module. We
define the \textit{cochain} \textit{complex} associated to the
module as an exact sequence:
$$
0\longrightarrow
C^0(\mathfrak{g},~\rm{A})\longrightarrow\cdots\longrightarrow
C^{q-1}(\mathfrak{g},~ \rm{A})\stackrel{\delta^{q-1}}\longrightarrow
C^{q}(\mathfrak{g},~\rm{A})\cdots.
$$

The spaces $C^q(\mathfrak{g},~\rm{A})$ are the spaces of super skew-symmetric $q$
linear mappings:
$$
C^0(\mathfrak{g},~\rm{A})=A,\quad C^q(\mathfrak{g},~\rm{A})=
\bigoplus_{q_0+q_1=q}\rm{Hom}(\bigwedge^{q_0}{\mathfrak g}_0\otimes S^{q_1}{\mathfrak g}_1,~\rm{A}).
$$
Elements of $C^q(\mathfrak{g},~\rm{A})$ are called \textit{cochains}.
The spaces $C^q(\mathfrak{g},~\rm{A})$ is $\mathbb{Z}_2$ graded:
$$
C^q(\mathfrak{g},~\rm{A})=C^q_0(\mathfrak{g},~\rm{A})+C^q_1(\mathfrak{g},~\rm{A}),
~~\text{with}~~C^q_p(\mathfrak{g},~\rm{A})= \bigoplus_{\begin{smallmatrix}q_0+q_1=
q\cr q_1+r=p~~ mod2\end{smallmatrix}}\rm{Hom}(\bigwedge^{q_0}\mathfrak{g}_0\otimes
S^{q_1}\mathfrak{g}_1,A_r).
$$
The linear mapping $\delta^q$ (or, briefly $\delta$) is called the
\textit{coboundary operator}. This operator is a generalization of
the usual Chevalley coboundary operator for Lie algebra to the
case of Lie superalgebra. Explicitly, it is defined as follows.
Take a cochain $c\in C^q(\mathfrak{g},~\rm{A})$, then for $q_0$,
$q_1$ with $q_0+q_1=q+1$, $\delta^qc$ is:

\begin{equation*}
\label{coboundary}
\begin{aligned}
\delta^qc(g_1,&~\ldots,~{g}_{q_0},~h_1,~\ldots,~h_{q_1})\cr
&=\sum_{1\leq s<t\leq q_0}(-1)^{s+t-1}c([g_s,g_t],g_1,\ldots,\hat
g_{s},\ldots, \hat g_{t},\ldots,g_{q_0},~h_1,\ldots,h_{q_1})\cr
&\hskip1cm+\sum_{s=1}^{q_0}\sum_{t=1}^{q_1}(-1)^{s-1}c(g_1,\ldots,\hat
g_{s}, \ldots,g_{q_0},~[g_s,h_t],h_1,\ldots,\hat
h_{t},\ldots,h_{q_1})\cr &\hskip1cm+\sum_{1\leq s<t\leq
q_1}c([h_s,h_t],g_1,\ldots,g_{q_0},~h_1,\ldots, \hat
h_{s},\ldots,\hat h_{t},\ldots,h_{q_1})\cr
&\hskip1cm+\sum_{s=1}^{q_0}(-1)^sg_sc(g_1,\ldots,\hat
g_s,\ldots,g_{q_0}, ~h_1,\ldots,h_{q_1})\cr
&\hskip1cm+(-1)^{q_{0}-1}\sum_{s=1}^{q_1}h_sc(g_1,\ldots,g_{q_0},~h_1,
\ldots,\hat h_{s},\ldots,h_{q_1}).
\end{aligned}
\end{equation*}
where $g_1,\ldots,g_{q_0}$ are in $\mathfrak{g}_0$ and $h_1,\ldots,
h_{q_1}$ in $\mathfrak{g}_1$.

The relation $\delta^q\circ\delta^{q-1}=0$ holds. The kernel of
$\delta^q$, denoted $Z^q(\mathfrak{g},~\rm{A})$, is the space of $q$
{\it cocycles}, among them, the elements in the range of
$\delta^{q-1}$ are called $q$ {\it coboundaries}. We note
$B^q(\mathfrak{g},~\rm{A})$ the space of $q$ coboundaries.

By definition, the $q^{th}$ cohomolgy space is the quotient space
$$
H^q (\mathfrak{g},~\rm{A})=Z^q(\mathfrak{g},~A)/B^q(\mathfrak{g},~\rm{A}).
$$

One can check that $\delta^q(C^q_p(\mathfrak{g},~\rm{A}))\subset C^{q+1}_p
(\mathfrak{g},~\rm{A})$ and then we get the following sequences
$$
0\longrightarrow C^0_p(\mathfrak{g},~\rm{A})\longrightarrow\cdots
\longrightarrow
C^{q-1}_p(\mathfrak{g},~\rm{A})\stackrel{\delta^{q-1}}
\longrightarrow C^q_p(\mathfrak{g},~\rm{A})\cdots,
$$
where $p=0$ or 1. The cohomology spaces are thus graded by
$$
H^q_p(\mathfrak{g},~\rm{A})=\rm{Ker}\delta^{q}|_{C^q_p(\mathfrak{g},
~\rm{A})}/\delta^{q-1}(C^{q-1}_p(\mathfrak{g},~\rm{A})).
$$

%------------------------------------------------------------------------------------------------
\subsection{The main theorem}
%-------------------------------------------------------------------------------------------------

The main result in this paper is the following:

\begin{thm}
\label{th1}

The cohomolgy spaces
$\mathrm{H}^1_p(\mathfrak{g},\mathfrak{D}_{\lambda,\mu})$ are
finite dimensional. An explicit description of these spaces is the
following:

\medskip
\noindent 1)The space
$\mathrm{H}^1_0(\mathfrak{osp}(1|2),\mathfrak{D}_{\lambda,\mu})$
is
\begin{equation}
\label{even cohom}
\mathrm{H}^1_0(\mathfrak{osp}(1|2),\mathfrak{D}_{\lambda,\mu})\simeq\left\{\begin{aligned}
\mathbb{R}&\hbox{ if }~~\lambda=\mu, \cr 0&\hbox{ otherwise.}
\end{aligned}
\right.
\end{equation}

A base for the space
$\mathrm{H}^1_{0}(\mathfrak{osp}(1|2),\mathfrak{D}_{\lambda,\lambda})$
is given by the cohomology class of the $1$-cocycle:
\begin{equation}
\label{even basis}
\Upsilon_{\lambda,\lambda}(X_F)=F'.
\end{equation}

\medskip
\noindent 2)The space
$\mathrm{H}^1_1(\mathfrak{osp}(1|2),\mathfrak{D}_{\lambda,\mu})$
is
\begin{equation}
\label{odd cohom}
\mathrm{H}^1_1(\mathfrak{osp}(1|2),\mathfrak{D}_{\lambda,\mu})\simeq\left\{\begin{aligned}
\mathbb{R}^2&\hbox{ if } ~~
\lambda=\frac{1-k}{2},~\mu=\frac{k}{2},\cr 0&\hbox{ otherwise. }
\end{aligned}
\right.
\end{equation}

A base for the space
$\mathrm{H}^1_{1}(\mathfrak{osp}(1|2),\mathfrak{D}_{\frac{1-k}{2},
\frac{k}{2}})$ is given by the cohomology classes of the
$1$-cocycles:
\begin{equation}
\label{odd basis}
\begin{aligned}
\Upsilon_{\frac{1-k}{2},\frac{k}{2}}(X_F)&=\overline{\eta}^2(F)\overline{\eta}^{2k-1},\cr
\widetilde{\Upsilon}_{\frac{1-k}{2},\frac{k}{2}}(X_F)&=(k-1)\eta^4(F)
\overline{\eta}^{2k-3}+\eta^3(F)\overline{\eta}^{2k-2}.
\end{aligned}
\end{equation}
\end{thm}

Note that the 1-cocycle $\widetilde{\Upsilon}_{\frac{1-k}{2},\frac{k}{2}}$
coincides with the 1-cocycle $\gamma_{2k-1}$ given by Gargoubi et
al. in \cite{gmo}. The proof of Theorem \ref{th1} will be the subject of subsection 3.4.

%-------------------------------------------------------------------------------------------------------------------
 \subsection{Relationship between $\mathrm{H}^1(\mathfrak{osp}(1|2),
 \mathfrak{D}_{\lambda,\mu})$ and $\mathrm{H}^1(\mathfrak{sl}(2),\mathcal{D}_{\lambda,\mu})$}

%-------------------------------------------------------------------------------------------------------------------------

Before proving the theorem \ref{th1} we present here some results illustrating the
analogy between the cohomlogy spaces in super and classical settings.

First, note that:
\begin{itemize}
\item[1)] As a $\mathfrak{sl}(2)$-module, we have $\mathfrak{F}_{\lambda}\simeq
\mathcal{F}_\lambda\oplus\Pi(\mathcal{F}_{\lambda+\frac{1}{2}})$ and
$\mathfrak{osp}(1|2)\simeq\mathfrak{sl}(2)\oplus\Pi(\mathfrak{h})$,
where $\mathfrak{h}$ is the subspace of $\mathcal{F}_{-\frac{1}{2}}$
spanned by $\{dx^{-\frac{1}{2}}, xdx^{-\frac{1}{2}}\}$ and $\Pi$ is
the change of parity.

\item[2)] As a $\mathfrak {sl}(2)$-module, we have for the homogeneous
components of $\mathfrak{D}_{\lambda,\mu}$:
\begin{equation*}
(\mathfrak{D}_{\lambda,\mu})_0\simeq\mathcal{D}_{\lambda,\mu}\oplus
\mathcal{D}_{\lambda+\frac{1}{2},\mu+\frac{1}{2}}\quad\hbox{and}\quad
(\mathfrak{D}_{\lambda,\mu})_1\simeq\Pi(\mathcal{D}_{\lambda+\frac{1}{2},\mu}
\oplus\mathcal{D}_{\lambda,\mu+\frac{1}{2}}).
\end{equation*}
\end{itemize}

\begin{proposition}\label{prop}
Any 1-cocycle $\Upsilon\in
Z^1(\mathfrak{{osp}}(1|2);\mathfrak{D}_{\lambda,\mu})$, is
decomposed into $(\Upsilon',\Upsilon'')$ in $Hom(\mathfrak{sl}(2);
\mathfrak{D}_{\lambda,\mu})\oplus
Hom(\mathfrak{h};\mathfrak{D}_{\lambda,\mu})$.
 $\Upsilon'$ and $\Upsilon''$ are solutions of the following equations:
\begin{align}
\label{coc1} &\Upsilon'([X_{g_1},X_{g_2}]) -\mathfrak{L}_{X_{g_1}}^{\lambda,\mu}\Upsilon'(X_{g_2})
+\mathfrak{L}_{X_{g_2}}^{\lambda,\mu} \Upsilon'(X_{g_1})= 0,\\
\label{coc2} &\Upsilon''([X_g,X_{h\theta}])-\mathfrak{L}_{X_{g}}^{\lambda,\mu}\Upsilon''(X_{h\theta})+
\mathfrak{L}_{X_{h\theta}}^{\lambda,\mu} \Upsilon'(X_{g})= 0,\\
\label{coc3} &\Upsilon'([X_{h_1\theta},X_{h_2\theta}])-
\mathfrak{L}_{X_{h_1\theta}}^{\lambda,\mu}\Upsilon''(X_{h_2\theta})-
\mathfrak{L}_{X_{h_2\theta}}^{\lambda,\mu}\Upsilon''(X_{h_1\theta})=0,
\end{align}
here, $g$, $g_1$, $g_2$ are polynomials in the variable $x$,
with degree at most 2, and $h$, $h_1$, $h_2$ are affine functions in the variable $x$.
\end{proposition}

\begin{proofname}.
The  equations (\ref{coc1}), (\ref{coc2}) and (\ref{coc3}) are
equivalent to the fact that $\Upsilon$ is a 1-cocycle. For any $X_F,\,X_G\in \mathrm{\mathfrak {osp}}(1|2)$,
\begin{equation*}
\delta \Upsilon(X_F,X_G):=\Upsilon([X_F,X_G])-
\mathfrak{L}_{X_F}^{\lambda,\mu}\Upsilon(X_G)+(-1)^{p(F)p(G)}\mathfrak{L}_{X_G}^{\lambda,\mu}\Upsilon(X_F)=0.
\end{equation*}
\end{proofname}
\hfill$\Box$

According to the $\mathbb{Z}_2$-grading, the even component $\Upsilon_0$ and the odd component $\Upsilon_1$ of any 1-cocycle $\Upsilon$ can be decomposed as $\Upsilon_0=(\Upsilon_{000},\Upsilon_{00\frac{1}{2}},\Upsilon_{110},\Upsilon_{11\frac{1}{2}})$ and  $\Upsilon_1=(\Upsilon_{010},\Upsilon_{01\frac{1}{2}},\Upsilon_{100},\Upsilon_{10\frac{1}{2}})$, where
\begin{align*}
\left\{\begin{array}{llll}
  \Upsilon_{000} :&\mathfrak{sl}(2)\hfill&\rightarrow&\mathcal{D}_{\lambda,\mu}, \\
  \Upsilon_{00\frac{1}{2}} :&\mathfrak{sl}(2)\hfill&\rightarrow&
  \mathcal{D}_{\lambda+\frac{1}{2},\mu+\frac{1}{2}}, \\
  \Upsilon_{110} :&\mathfrak{h}\hfill&\rightarrow&
  \mathcal{D}_{\lambda,\mu+\frac{1}{2}}, \\
  \Upsilon_{11\frac{1}{2}} :&\mathfrak{h}\hfill&\rightarrow&
  \mathcal{D}_{\lambda+\frac{1}{2},\mu}
 \end{array}\right.
~~\hbox{ and }~~\left\{\begin{array}{llll}
  \Upsilon_{010} :&\mathfrak{sl}(2)\hfill&\rightarrow&
  \mathcal{D}_{\lambda,\mu+\frac{1}{2}}, \\
  \Upsilon_{01\frac{1}{2}} :&\mathfrak{sl}(2)\hfill&
  \rightarrow&\mathcal{D}_{\lambda+\frac{1}{2},\mu} \\
  \Upsilon_{100} :&\mathfrak{h}\hfill&\rightarrow&
  \mathcal{D}_{\lambda,\mu}, \\
  \Upsilon_{10\frac{1}{2}}:&\mathfrak{h}\hfill&\rightarrow&
  \mathcal{D}_{\lambda+\frac{1}{2},\mu+\frac{1}{2}}.
 \end{array}\right.
\end{align*}
The decomposition $\Upsilon=(\Upsilon',\Upsilon'')$ given in proposition \ref{prop}
corresponds to
$$
\Upsilon'=(\Upsilon_{000},\Upsilon_{00\frac{1}{2}},\Upsilon_{010},
\Upsilon_{01\frac{1}{2}})\quad\text{and}\quad \Upsilon''=(\Upsilon_{110},
\Upsilon_{11\frac{1}{2}},\Upsilon_{100},\Upsilon_{10\frac{1}{2}}).
$$

By considering the equation (\ref{coc1}), we can see the components
$\Upsilon_{000}$, $\Upsilon_{00\frac{1}{2}} $, $\Upsilon_{010}$ and
$\Upsilon_{01\frac{1}{2}}$ as $1$-cocycles on $\mathfrak{sl}(2)$
with coefficients respectively in $\mathcal{D}_{\lambda,\mu}$,
$\mathcal{D}_{\lambda+\frac{1}{2},\mu+\frac{1}{2}}$,
$\mathcal{D}_{\lambda,\mu+\frac{1}{2}},$ and
$\mathcal{D}_{\lambda+\frac{1}{2},\mu}$.

The first cohomology space
$\mathrm{H}^1(\mathfrak{sl}(2);\mathcal{D}_{\lambda,\mu})$ was
computed by Gargoubi and Lecomte \cite{g,lec}. The result is the
following:
\begin{equation}
\mathrm{H}^1(\mathfrak{sl}(2);\mathcal{D}_{\lambda,\mu})\simeq\left\{
\begin{array}{llll}
\mathbb{R}&\hbox{if}&\lambda=\mu&\\
\mathbb{R}^2&\hbox{if}&(\lambda,\mu)=(\frac{1-k}{2},\frac{1+k}{2})&
\text{where}\quad k\in\mathbb{N}\setminus\{0\}\\
0&&\hbox{otherwise}.&
\end{array}
\right.
\end{equation}
The space
$\mathrm{H}^1(\mathfrak{sl}(2);\mathcal{D}_{\lambda,\lambda})$ is
generated by the cohomology class of the 1-cocycle
\begin{equation}\label{slcocy0}
C'_\lambda(F\frac{d}{dx})(f{dx}^{\lambda})=F'f{dx}^{\lambda}.
\end{equation}
For $k\in\mathbb{N}\setminus\{0\}$, the space
$\mathrm{H}^1(\mathfrak{sl}(2);
\mathcal{D}_{\frac{1-k}{2},\frac{1+k}{2}})$ is generated by the
cohomology classes of the 1-cocycles, $C_k$ and $\tilde{C}_k$
defined by
\begin{align}
\label{slcocy1}
C_k(F\frac{d}{dx})(f{dx}^{\frac{1-k}{2}})=F'f^{(k)}{dx}^{\frac{1+k}{2}}\quad\text{ and }
\quad\widetilde{C}_k(F\frac{d}{dx})(f{dx}^{\frac{1-k}{2}})=
F''f^{(k-1)}dx^{\frac{1+k}{2}}.
\end{align}

We shall need the following description of $\mathfrak{sl}(2)$ invariant mappings.
\begin{lem}\label{inv}
Let
\begin{equation*}
A:\mathfrak{h}\times\mathcal{F}_\lambda\rightarrow\mathcal{F}_\mu,
\qquad(hdx^{-\frac{1}{2}},fdx^\lambda)\mapsto A(h,f)dx^{\mu}
\end{equation*}
be a bilinear differential operator. If $A$ is
$\mathfrak{sl}(2)$-invariant then
\begin{equation*}
\mu=\lambda-\frac{1}{2}+k,\quad\text{where}\quad k\in\mathbb{N}
\end{equation*}
and the following relation holds
\begin{equation*}
A_k(h,f)=a_k(h f^{(k)} +k(2\lambda+k-1)h'f^{(k-1)}),\quad\text{where}
\quad k(k-1)(2\lambda+k-1)(2\lambda+k-2)a_k=0.
\end{equation*}
\end{lem}

\begin{proofname}. A straightforward computation.
\end{proofname}

\hfill$\Box$

Now, let us study the relationship between these 1-cocycles and
their analogues in the super setting. We know that any element
$\Upsilon\in Z^1(\mathfrak{osp}(1|2), \mathfrak{D}_{\lambda,\mu})$
is decomposed into $\Upsilon=\Upsilon'+\Upsilon''$ where
$\Upsilon'\in
Hom\left(\mathfrak{sl}(2),\mathfrak{D}_{\lambda,\mu}\right)$ and
$\Upsilon''\in
Hom\left(\mathfrak{h},\mathfrak{D}_{\lambda,\mu}\right)$. The
following lemma shows the close relationship between the cohomolgy
spaces $H^1(\mathfrak{osp}(1|2),\mathfrak{D}_{\lambda,\mu})$ and
$H^1(\mathfrak{sl}(2),\mathcal{D}_{\lambda,\mu})$.

\begin{lem}\label{sa}
The 1-cocycle $\Upsilon$ is a coboundary for $\mathfrak{osp}(1|2)$ if and only if
$\Upsilon'$ is a coboundary for $\mathfrak{sl}(2)$.
\end{lem}

\begin{proofname}. It is easy to see that if $\Upsilon$ is a coboundary for
$\mathfrak{osp}(1|2)$ then  $\Upsilon'$ is a coboundary over $\mathfrak{sl}(2)$.
Now, assume that $\Upsilon'$ is a coboundary for $\mathfrak{sl}(2)$, that is,
there exists $\widetilde{A}\in\mathfrak{D}_{\lambda,\mu}$ such that for all $g$
polynomial in the variable $x$ with degree at most 2
\begin{equation*}
\Upsilon'(X_g)=\mathfrak{L}_{X_g}^{\lambda,\mu}\widetilde{A}.
\end{equation*}
By replacing $\Upsilon$ by $\Upsilon-\delta\widetilde{A}$, we can suppose that
$\Upsilon'=0$. But, in this case, the map $\Upsilon''$ must satisfy, for all
$h$, $h_1$, $h_2$ polynomial with degree 0 or 1 and $g$ polynomial with degree 0,1 or 2,
the following equations
\begin{align}
 \label{sltr1} &\mathfrak{L}^{\lambda,\mu}_{X_g}\Upsilon''(X_{h\theta})-
 \Upsilon''([X_g,X_{h\theta}])=0, \\
 \label{sltr2} &\mathfrak{L}^{\lambda,\mu}_{X_{h_1\theta}}\Upsilon''(X_{h_2\theta}) +
 \mathfrak{L}^{\lambda,\mu}_{X_{h_2\theta}}      \Upsilon''(X_{h_1\theta})=0.
\end{align}
1) If $\Upsilon$ is an even 1-cocycle then $\Upsilon''$ is
decomposed into $\Upsilon''_{00}:\mathfrak{h}\otimes
\mathcal{F}_{\lambda+\frac{1}{2}} \rightarrow\mathcal{F}_\mu$ and
$\Upsilon''_{01}:\mathfrak{h}\otimes
\mathcal{F}_\lambda\rightarrow\mathcal{F}_{\mu+\frac{1}{2}}$. The
equation (\ref{sltr1}) tell us that $\Upsilon''_{00}$ and
$\Upsilon''_{01}$ are $\mathfrak{sl}(2)$ invariant bilinear maps.
Therefore, the expressions of $\Upsilon''_{00}$ and
$\Upsilon''_{01}$ are given by Lemma \ref{inv}. So, we must have
$\mu=\lambda+k=(\lambda+\frac{1}{2})-\frac{1}{2}+k$ (and then
$\mu+\frac{1}{2}=\lambda-\frac{1}{2}+k+1$). More precisely, using
the equation (\ref{sltr2}), we get up to a factor:
\begin{equation*}
\Upsilon=\left\{\begin{aligned}
&0\hskip 1.4cm\text{ if }~k(k-1)(2\lambda+k)(2\lambda+k-1)\neq0~~\text{ or }~k=1
\text{ and }\lambda\notin\{0,\,-\frac{1}{2}\},\\
&\delta(\theta\partial_\theta\partial_x^k)\hskip0.3cm\text{ if }~(\lambda,\mu)=
(\frac{-k}{2},\frac{k}{2}),\\
&\delta (\partial_x^k-\theta\partial_\theta\partial_x^k)~~\text{ if }~(\lambda,\mu)=
(\frac{1-k}{2},\frac{1+k}{2})~~\text{ or }~\lambda=\mu.
\end{aligned}\right.
\end{equation*}
2) If $\Upsilon$ is an odd 1-cocycle then $\Upsilon''$ is decomposed into
$\Upsilon''_{00}:\mathfrak{h}\otimes \mathcal{F}_{\lambda}\rightarrow \mathcal{F}_\mu$ and
$\Upsilon''_{01}:\mathfrak{h}\otimes\mathcal{F}_{\lambda+\frac{1}{2}}
\rightarrow\mathcal{F}_{\mu+\frac{1}{2}}$.
As in the previous case, the expressions of $\Upsilon''_{00}$ and $\Upsilon''_{01}$
are given by Lemma \ref{inv}. So, we must have $\mu=\lambda-\frac{1}{2}+k$
(and then $\mu+\frac{1}{2}=(\lambda+\frac{1}{2})-\frac{1}{2}+k$.)
More precisely, using the equation (\ref{sltr2}), we get:
\begin{equation*}
\Upsilon=\left\{\begin{array}{lll}
0&\text{ if }& k(k-1)(2\lambda+k-1)\neq0\\
\delta(\theta)&\text{ if }&\mu=\lambda-\frac{1}{2},\\
\delta(\partial_\theta)&\text{ if }&\mu=\lambda+\frac{1}{2},\\
\delta(\theta\partial_x^k)&\text{ if }&(\lambda,\mu)=(\frac{1-k}{2},\frac{k}{2}).
\end{array}\right.
\end{equation*}
\end{proofname}
\hfill$\Box$

Now, the space $Z^1(\mathfrak{osp}(1|2),\mathfrak{D}_{\lambda,\mu})$ of
1-cocycles is $\mathbb{Z}_2$-graded:
\begin{equation}\label{graded}
Z^1(\mathfrak{osp}(1|2),\mathfrak{D}_{\lambda,\mu})=
Z^1(\mathfrak{osp}(1|2),\mathfrak{D}_{\lambda,\mu})_0\oplus
Z^1(\mathfrak{osp}(1|2), \mathfrak{D}_{\lambda,\mu})_1.
\end{equation}
Therefore, any element $\Upsilon\in Z^1(\mathfrak{osp}(1|2),\mathfrak{D}_{\lambda,\mu})$
is decomposed into an even part $\Upsilon_0$ and odd part $\Upsilon_1$. Each of
$\Upsilon_0$ and $\Upsilon_1$ is decomposed into two components:
$\Upsilon_0=(\Upsilon_{00},\Upsilon_{11})$ and
$\Upsilon_1=(\Upsilon_{01},\Upsilon_{10})$, where
\begin{align*}
\left\{\begin{array}{ll}
  \Upsilon_{00}:\mathfrak{sl}(2)&\rightarrow(\mathfrak{D}_{\lambda,\mu})_0, \\
  \Upsilon_{11}:\mathfrak{h}&\rightarrow(\mathfrak{D}_{\lambda,\mu})_1,
 \end{array}\right.
\hbox{ and }
\left\{\begin{array}{ll}
  \Upsilon_{01}:\mathfrak{sl}(2)&\rightarrow(\mathfrak{D}_{\lambda,\mu})_1, \\
  \Upsilon_{10}:\mathfrak{h}&\rightarrow(\mathfrak{D}_{\lambda,\mu})_0.
 \end{array}\right.
\end{align*}
The components $\Upsilon_{11}$ and $\Upsilon_{10}$ of $\Upsilon_0$ and
$\Upsilon_1$ are also decomposed as follows:
$\Upsilon_{11}=\Upsilon_{110}+\Upsilon_{11\frac{1}{2}}$ and
$\Upsilon_{10}=\Upsilon_{100}+\Upsilon_{10\frac{1}{2}}$, where
$\Upsilon_{110}\in Hom\left(\mathfrak{h},\mathcal{D}_{\lambda,\mu+\frac{1}{2}}\right)$,
$\Upsilon_{11\frac{1}{2}}\in Hom\left(\mathfrak{h}, \mathcal{D}_{\lambda+\frac{1}{2},\mu}\right)$,
$\Upsilon_{100}\in Hom\left(\mathfrak{h},\mathcal{D}_{\lambda,\mu}\right)$,
$\Upsilon_{10\frac{1}{2}}\in
Hom\left(\mathfrak{h},\mathcal{D}_{\lambda+\frac{1}{2},\mu+\frac{1}{2}}\right)$.

As in \cite{bbbbk}, the following lemma gives the general form of each of
$\Upsilon_{110}$ and $\Upsilon_{11\frac{1}{2}}$.

\begin{lem}\label{sd}
Up to a coboundary, the maps $\Upsilon_{110}$, $\Upsilon_{11\frac{1}{2}}$,
$\Upsilon_{100}$ and $\Upsilon_{10\frac{1}{2}}$ are given by
\begin{equation*}
\begin{aligned}
\Upsilon_{110}(X_{h\theta})&=a_0h\theta\partial^k_x+a_1h'\theta\partial^{k-1}_x
~~\hbox{and}~~\Upsilon_{11\frac{1}{2}}(X_{h\theta})&=b_0
h\partial_{\theta}\partial^k_x+b_1h'\partial_{\theta}\partial^{k-1}_x,\\
\Upsilon_{110}(X_{h\theta})&=c_0h\theta\partial^k_x+c_1h'\theta\partial^{k-1}_x
~~\hbox{and}~~\Upsilon_{11\frac{1}{2}}(X_{h\theta})&=d_0
h\partial_{\theta}\partial^k_x+d_1h'\partial_{\theta}\partial^{k-1}_x,
\end{aligned}
\end{equation*}
where the coefficients $a_i$, $b_i$, $c_i$, and $d_i$ are constants.
\end{lem}

\begin{proofname}.
The coefficients $a_i$, $b_i$, $c_i$, and $d_i$ a priori are some functions of
$x$, but we shall now prove $\partial_xa_i=\partial_xb_i=0$ (and similarly
$\partial_xc_i=\partial_xd_i=0$). To do that, we shall simply show that
$\mathfrak{L}^{\lambda,\mu}_{\partial_x}(\Upsilon_{11})=0$.

First, for all $h$ polynomial with degree 0 or 1, we have
\begin{equation}\label{partial}
(\mathfrak{L}^{\lambda,\mu}_{\partial_x}\Upsilon_{11})(X_{h\theta})=
\mathfrak{L}^{\lambda,\mu}_{\partial_x}(\Upsilon_{11}(X_{h\theta}))-
\Upsilon_{11}([\partial_x,X_{h\theta}]).
\end{equation}
On the other hand, from Lemma \ref{sa}, it follows that, up to a coboundary,
$\Upsilon_{00}$ is a linear combination of some 1-cocycles for $\mathrm{sl}(2)$
given by (\ref{slcocy0}) and (\ref{slcocy1}). So, we have $\Upsilon_{00}(\partial_x)=0$
and then
\begin{equation*}
\mathfrak{L}^{\lambda,\mu}_{X_{h\theta}}(\Upsilon_{00}(\partial_x))=0.
\end{equation*}
Therefore, the equation (\ref{partial}) becomes, for all $h$,
\begin{equation}\label{partial1}
-(\mathfrak{L}^{\lambda,\mu}_{\partial_x}\Upsilon_{11})(X_{h\theta})=
\Upsilon_{11}([\partial_x,X_{h\theta}])-
\mathfrak{L}^{\lambda,\mu}_{\partial_x}(\Upsilon_{11}(X_{h\theta}))+
\mathfrak{L}^{\lambda,\mu}_{X_{h\theta}}(\Upsilon_{00}(\partial_x)).
\end{equation}

The right-hand side of (\ref{partial1}) is nothing but
$\delta\Upsilon_0(\partial_x,X_{h\theta})$. But, $\Upsilon_0$ is a
1-cocycle, then
$\mathfrak{L}^{\lambda,\mu}_{\partial_x}(\Upsilon_{11})=0.$ Lemma
\ref{sd} is proved.
\end{proofname}

\hfill$\Box$

%-----------------------------------------------------------------------------------------------------

\subsection{Proof of Theorem \ref{th1}}
%---------------------------------------------------------------------------------------------------------

The first cohomology space
$\mathrm{H}^1(\mathfrak{osp}(1|2);\mathfrak{D}_{\lambda,\mu})$
inherits the $\mathbb{Z}_2$-grading from (\ref{graded}) and is
decomposed into odd and an even subspaces:
$$
\mathrm{H}^1(\mathfrak {osp}(1|2);\mathfrak{D}_{\lambda,\mu})=
\mathrm{H}^1_0(\mathfrak{osp}(1|2);\mathfrak{D}_{\lambda,\mu})\oplus
\mathrm{H}^1_1(\mathfrak{osp}(1|2);\mathfrak{D}_{\lambda,\mu}).
$$
We compute each part separetly.

1) Let $\Upsilon_0$ be a non trivial even 1-cocycle for $\mathfrak{osp}(1|2)$ in
$\mathfrak{D}_{\lambda,\mu}$. According to the $\mathbb{Z}_2$-grading, $\Upsilon_0$
should retain the following general form: $\Upsilon_0=\Upsilon_{000}+$
$\Upsilon_{00\frac{1}{2}}+\Upsilon_{110}+\Upsilon_{11\frac{1}{2}}$
such that
\begin{equation}\label{decomp1}
\left\{\begin{array}{lllll}
\Upsilon_{000}&:& \mathfrak{sl}(2)&\rightarrow&\mathcal{D}_{\lambda,\mu}, \\
\Upsilon_{00\frac{1}{2}}&:& \mathfrak{sl}(2)&\rightarrow&\mathcal{D}_{\lambda\frac{1}{2}f,\mu+\frac{1}{2}}, \\
\Upsilon_{110}&:&~~~~\mathfrak{h} &\rightarrow&\mathcal{D}_{\lambda,\mu+\frac{1}{2}}, \\
\Upsilon_{11\frac{1}{2}}&:&~~~~ \mathfrak{h} &\rightarrow&
 \mathcal{D}_{\lambda+\frac{1}{2},\mu}.
\end{array}\right.\end{equation}
Then, by using Lemma \ref{sa}, we deduce that, up to coboundary, $\Upsilon_{000}$
and $\Upsilon_{00\frac{1}{2}}$ can be expressed in terms of $C'_{\lambda}$, $C_k$
and $\widetilde{C}_k$ where $\lambda\in\mathbb{R}$ and $k\in\mathbb{N}\setminus\{0\}$.
We thus consider three cases:
\begin{itemize}
  \item [i)] $\lambda=\mu$,  $\Upsilon_{000}=\alpha C'_{\lambda}$, and
  $\Upsilon_{00\frac{1}{2}}=\beta C'_{\lambda\frac{1}{2}f}$.
  \item [ii)] $(\lambda,\mu)=(\frac{1-k}{2},\frac{1+k}{2})$,
  $\Upsilon_{000}=\alpha_1 C_{k}+\alpha_2\widetilde{C}_k$, and
  $\Upsilon_{00\frac{1}{2}}=0$.
  \item [iii)] $(\lambda,\mu)=(\frac{-k}{2},\frac{k}{2})$, $\Upsilon_{000}=0$
  and $\Upsilon_{00\frac{1}{2}}=\alpha_1 C_{k}+\alpha_2\widetilde{C}_k$.
\end{itemize}
Put $\Upsilon'=\Upsilon_{000}+\Upsilon_{00\frac{1}{2}}$ and
$\Upsilon''=\Upsilon_{110}+\Upsilon_{11\frac{1}{2}}$. In each case,
the 1-cocycle $\Upsilon_0$ must satisfy
\begin{equation}\label{decomp5}
\left\{\begin{array}{lllll}
\Upsilon''[X_g,X_{\theta h}]&=&\mathfrak{L}_{X_g}^{\lambda,\mu}\Upsilon''(X_{\theta h})-
\mathfrak{L}_{X_{\theta h}}^{\lambda,\mu}\Upsilon'(X_g), \\
\Upsilon'[X_{\theta h_1},X_{\theta h_2}]&=&\mathfrak{L}_{X_{\theta h_1}}^{\lambda,\mu}
\Upsilon''(X_{\theta h_2})+\mathfrak{L}_{X_{\theta h_2}}^{\lambda,\mu}\Upsilon''(X_{\theta h_1}),
\end{array}\right.\end{equation}
where $h$, $h_1$, and $h_2$ are polynomials of degree 0 or 1, $g$
polynomial of degree 0, 1 or 2.

Now, thanks to Lemma \ref{sd}, we can write
\begin{equation*}
\Upsilon_{110}(X_{h\theta})=a_0h\theta\partial^k_x~+a_1h'\theta\partial^{k-1}_x~~\hbox{and}
~~\Upsilon_{11\frac{1}{2}}(v_{h\theta})=
b_0h\partial_{\theta}\partial^k_x+b_1h'\partial_{\theta}\partial^{k-1}_x.
\end{equation*}
Let us now solve the equations (\ref{decomp5}). We obtain
$\lambda=\mu$ and $\Upsilon_{\lambda,\lambda}(X_F)=F'$.
This completes the proof of part 1).

\vskip.3cm
2) Consider a non trivial odd 1-cocycle $\Upsilon_1$ for
$\mathfrak{osp}(1|2)$ in $\mathfrak{D}_{\lambda,\mu}$ and
its decomposition $\Upsilon_1=\Upsilon_{010}+\Upsilon_{01\frac{1}{2}}+
\Upsilon_{100}+\Upsilon_{10\frac{1}{2}}$, where
\begin{equation}\label{decomp3}
\left\{\begin{array}{lllll}
\Upsilon_{010}&:&\mathfrak{sl}(2)&\rightarrow &\mathcal{D}_{\lambda,\mu+\frac{1}{2}}, \\
\Upsilon_{01\frac{1}{2}}&:&\mathfrak{sl}(2)&\rightarrow &
\mathcal{D}_{\lambda+\frac{1}{2},\mu},\\
\Upsilon_{100}&:& ~~~~\mathfrak{h}&\rightarrow &\mathcal{D}_{\lambda,\mu}, \\
\Upsilon_{10\frac{1}{2}}&:&~~~~\mathfrak{h}&\rightarrow &
\mathcal{D}_{\lambda+\frac{1}{2},\mu+\frac{1}{2}}.
\end{array}\right.
\end{equation}
We must have $(\lambda,\mu)=(\frac{1-k}{2},\frac{k}{2})$ with
$k\in\mathbb{N}\setminus\{0\}$. Moreover
$\Upsilon_1=\Upsilon_{000}+\Upsilon_{00\frac{1}{2}} +
\Upsilon_{110}+\Upsilon_{11\frac{1}{2}}$ is a 1-cocycle for
$\mathcal{K}(1)$ if and only if
\begin{equation}\label{decomp6}
\left\{\begin{array}{lllll}
\Upsilon_{010}&=&\alpha_1 C_{k}+\alpha_2\widetilde{C}_k\\
\Upsilon_{01\frac{1}{2}}&=&\beta_1 C_{k-1}+\beta_2\widetilde{C}_{k-1}\\
\Upsilon''[X_g,\,X_{\theta h}]&=&\mathfrak{L}_{X_g}^{\lambda,\mu}\Upsilon''
(X_{\theta h})-\mathfrak{L}_{X_{\theta h}}^{\lambda,\mu}\Upsilon'(X_g), \\
\Upsilon'[X_{\theta h_1},\,X_{\theta h_2}]&=&\mathfrak{L}_{X_{\theta h_1}}^{\lambda,\mu}
\Upsilon''(X_{\theta h_2}) +\mathfrak{L}_{X_{\theta
h_2}}^{\lambda,\mu}\Upsilon''(X_{\theta h_1}),
\end{array}\right.\end{equation}
where $\Upsilon'=\Upsilon_{010}+\Upsilon_{01\frac{1}{2}}$ and
$\Upsilon''=\Upsilon_{100}+\Upsilon_{10\frac{1}{2}}$.

As above, we then can write
\begin{equation*}
\Upsilon_{100}(X_{h\theta})=a_0h\theta\partial^k_x~+a_1h'\theta\partial^{k-1}_x
~~\hbox{and}~~\Upsilon_{10\frac{1}{2}}(v_{h\theta})=b_0h
\partial_{\theta}\partial^k_x+b_1h'\partial_{\theta}\partial^{k-1}_x.
\end{equation*}
According to Lemma \ref{sa}, the map $\Upsilon_1$ is a non trivial
1-cocycle if and only if at least one of the maps $\Upsilon_{010}$ and
$\Upsilon_{01\frac{1}{2}}$ is a non trivial 1-cocycle for $\mathfrak{sl}(2)$,
that means $(\alpha_1,\alpha_2,\beta_1,\beta_2)\neq(0,0,0,0)$. Let us determine
the linear maps $\Upsilon_{100}$ and $\Upsilon_{10\frac{1}{2}}$. Up to factor, we get:
\begin{equation*}
\Upsilon_1=\alpha_1\Upsilon_{\frac{1-k}{2},\frac{k}{2}}+
\alpha_2\widetilde{\Upsilon}_{\frac{1-k}{2},\frac{k}{2}}+a_0\delta(2\theta\partial_x^k).
\end{equation*}
Thus, the cohomology classes of $\Upsilon_{\frac{1-k}{2},\frac{k}{2}}$ and
$\widetilde{\Upsilon}_{\frac{1-k}{2},\frac{k}{2}}$ generate
$H^1_{1}(\mathfrak{osp}(1|2),\mathfrak{D}_{\frac{1-k}{2},\frac{k}{2}})$.
The proof is now complete.

%---------------------------------------------------------------------------------------------------------------------

%--------------------------------------------------------------------------------------------------------------


\begin{thebibliography}{99}

%---------------------------------------------------------------------------------------------------------------------
\bibitem{bbbbk}
I. Basdouri, M. Ben Ammar, N. Ben Fraj, M. Boujelbene and K. Kammoun {\it Cohomology of the Lie Superalgebra of Contact Vector Fields on
$\mathbb{R}^{1|1} $ and Deformations of the Superspace of Symbols}, ~math.RT/0702645.

\bibitem{Fu}
Fuchs  D B, {\it Cohomology of infinite-dimensional Lie algebras}, Plenum Publ. New York, 1986.

\bibitem{g}
H. Gargoubi, {\it Sur la g\'eom\'etrie de l'espace des op\'erateurs
diff\'erentiels lin´eaires sur $\mathbb{R}$}, Bull. Soc. Roy. Sci.
Li\`ege. Vol. 69, 1, 2000, 21–47.


\bibitem{gmo} H. Gargoubi, N. Mellouli and V. Ovsienko {\it Differential Operators on Supercircle: Conformally Equivariant Quantization and Symbol Calculus}, Letters in Mathematical Physics (2007) {\bf 79}: 51–65.

\bibitem{lec}
P. B. A. Lecomte, {\it On the cohomology of $\mathfrak{sl}(n+1;\mathbb{R})$ acting on differential operators and $\mathfrak{sl}(n+1;\mathbb{R})$-equivariant symbols}, Indag. Math. NS. 11 (1), (2000), 95 114.

\bibitem{nr}
A. Nijenuis, R. W. Richardson Jr., {\it Deformations of homomorphisms of Lie groups and Lie algebras}, Bull. Amer. Math. Soc. {\bf 73} (1967), 175--179.

\end{thebibliography}
\end{document}